\documentstyle[12pt]{article}
\newtheorem{thm}{Theorem}
\newtheorem{lm}{Lemma}

\newcommand{\be}{\begin{equation}}
\newcommand{\ee}{\end{equation}}
\newcommand{\bea}{\begin{eqnarray}}
\newcommand{\eea}{\end{eqnarray}}
\newcommand{\ba}{\begin{array}}
\newcommand{\ea}{\end{array}}
\newcommand{\beas}{\begin{eqnarray*}}
\newcommand{\eeas}{\end{eqnarray*}}
\newcommand{\ZZ}{{\bf Z}}
\newcommand{\RR}{{\bf R}}

\newcommand{\QQ}{{\bf Q}}
\newcommand{\DD}{{\cal D}}
\newcommand{\FF}{{\cal F}}
\newcommand{\nn}{\nonumber}

\def\vep{\varepsilon}

\def\xit{\widetilde{\xi}}
\begin{document}

\large
%\Large
\title{A variational coupling for a  \\ totally asymmetric
 exclusion process \\ with long jumps but no passing }

\normalsize

\author{Timo Sepp\"al\"ainen \\
Department of Mathematics, Iowa State University\\
Ames, Iowa 50011, USA\\
seppalai@iastate.edu\\
}

\maketitle

\begin{abstract} 
We prove a weak law of large numbers for a tagged particle in
a totally asymmetric exclusion process on the one-dimensional
lattice. The particles are allowed to take long jumps but not
pass each other. The object of the paper is to illustrate
a special technique for proving such theorems. The method uses
a coupling that mimics the Hopf-Lax formula from the theory of
viscosity solutions of Hamilton-Jacobi equations. 
\end{abstract} 

\vspace{1.cm}


\noindent
Research supported in part by National Science Foundation
grant DMS-9801085.

\newpage
\section{Introduction}

The purpose of this paper is to explain a technique for
deriving scaling limits that applies to
 certain totally asymmetric particle
processes. The example we use to illustrate the method is a
one-dimensional 
exclusion process whose particles are allowed to take 
arbitrarily long jumps, as long as particles do not pass
each other. Thus the ordering of particles is preserved,
and the process models traffic
on a single-lane highway. To the author's knowledge, 
the hydrodynamics of this particular process have not
been addressed before. 

The key point of the technique is this: There exists 
a coupling that expresses a process with
general initial conditions in terms of a family of processes
that have a simple structure. The evolution 
from a general initial configuration is constructed from
the simple cases by a variational formula, hence the
term ``variational coupling''. The benefit of the scheme is that
the limit behavior of the simple processes is tractable,
and through the coupling  one obtains results for the 
general process. 

This paper is intended as an introduction to the technique. 
Fairly complete proofs are provided, and we have not strived
to obtain the best possible result for the process under
study. Some  extensions of the theorem of the paper
and shortcomings of the method are addressed in
Section \ref{comments}. 
A review of limit theorems for tagged particles and 
references to literature can be found in \cite{fe}.

{\it Notational remarks.} For a real number $x$, $[x]$
denotes the largest integer less than or equal to $x$. 
$\ZZ_-$ is the set of nonpositive integers
$\{\ldots,-3,-2,-1,0\}$.

\section{Description and construction of the process}
\label{constr}

The process consists of particles that move on the one-dimensional lattice
$\ZZ$. At most one particle is permitted on each site. (This is
the meaning of the name ``exclusion process''.) We label the 
particles by integers $i\in\ZZ$. The state of the process
is a sequence of integers $\sigma=(\sigma(i):i\in\ZZ)$,
 where $\sigma(i)$ is the location of particle $i$. 
We shall also use ``$\sigma(i)$'' as the name of the particle,
because typically we need to 
discuss  simultaneously
 several different processes denoted by different Greek letters. 
The dynamics will be such that
 the ordering of particles is preserved. So we can 
assume that 
\begin{equation}
 \sigma(i)+1\le \sigma(i+1) \mbox{ for all $i\in\ZZ$.}
\label{e1}
\end{equation}
This condition contains both the exclusion rule  and the 
no-passing rule. 

It is also convenient to allow the values $\pm \infty$ for
$\sigma(i)$. Assumption  (\ref{e1}) remains valid.
 If $\sigma(i)<\infty=\sigma(i+1)$,
then $\sigma(i)$ is the rightmost particle on the lattice
$\ZZ$ because $\sigma(j)=\infty$ for all $j>i$. Similarly
there is a leftmost particle if some $\sigma(j)=-\infty$. 

To summarize, the state space of the process is the set
of $\ZZ\cup\{\pm \infty\}$-valued sequences $\sigma$
that satisfy  (\ref{e1}). 

For the dynamics, assume given a sequence $\beta_k$, $k\ge 1$,
of nonnegative real numbers. These are the rates. 
For each $k$ and each $i$, 
 particle $\sigma(i)$ independently
 attempts to jump $k$ steps to the 
right at exponential rate $\beta_k$. If there are $k$ empty lattice spaces
in front of $\sigma(i)$, this jump can be executed. If there
are not, in other words, if $\sigma(i+1)\le \sigma(i)+k$, then
$\sigma(i)$ jumps as much as it can, to the location 
$\sigma(i+1)-1$ right behind the next particle. The jumps
occur independently for all particles $\sigma(i)$ and all
jump sizes $k$.

 For nonnegative integers $m$, define  
\be
B_m=\sum_{k=1}^\infty k^m\beta_k\,.
\label{defBm} 
\ee
We assume 
\be
B_2<\infty.
\label{assB2} 
\ee
Then also $B_0<\infty$,
 which guarantees that there are finitely many jump attempts
in finite time intervals almost surely.  $B_2<\infty$ implies that 
 the distance traveled by any 
particle during a finite time interval has a finite second moment. 

To rigorously construct the process, let 
$\{\DD_i: i\in\ZZ\}$ be an i.i.d.\  collection
of homogeneous rate $B_0$ Poisson point processes on $[0,\infty)$.
In other words, a time interval $[t,t+h)$ has a Poisson($B_0h$)-distributed
random number of points (or {\it epochs}) from the process
 $\DD_i$. For each process $\DD_i$, label the points 
independently by integers
$k\ge 1$, with probability $\beta_k/B_0$ for label $k$. 
Let $\DD^k_i$ denote the point process of $k$-epochs
(that is, the points with label $k$) from 
process $\DD_i$. It follows that 
$\{\DD^k_i: i\in\ZZ,\,k\ge 1\}$ is a collection
of mutually  independent 
homogeneous Poisson point processes on $[0,\infty)$, and 
the rate of $\DD^k_i$ is $\beta_k$.
 If $\beta_k=0$ the point
process $\DD^k_i$ is empty for all $i$.

 The idea of the construction is that 
particle $\sigma(i)$ reads its jump commands from the 
Poisson process $\DD_i$, and attempts 
to execute a $k$-step jump at each epoch of $\DD^k_i$. 
Let $\sigma(i,t)$ denote the location of particle
$\sigma(i)$ at time $t\ge 0$, and 
$\sigma(t)=(\sigma(i,t):i\in\ZZ)$ the entire state  at time $t$.
 Assume that a
deterministic  initial configuration 
$\sigma(0)=(\sigma(i,0):i\in\ZZ)$ has been specified. 
Informally, the rule of evolution is this:

\paragraph{Jump rule.}
Suppose $\tau$ is an epoch of $\DD^k_i$. Assume that
the state $\sigma(t)$ has been defined for
$t<\tau$. At time $\tau$ set
\be
\sigma(i,\tau)=\min\left\{ \sigma(i,\tau-)+k\,,\,\sigma(i+1,\tau-)-1
\right\}\,.
\label{jump2}
\ee
After $\tau$, $\sigma(i)$ stays constant until the next epoch 
in $\DD_i$. The same scheme happens simultaneously and 
independently for all particles $\sigma(i)$. 

\hbox{}

Since there are infinitely many particles,
infinitely many jumps are attempted 
 almost surely in any 
time interval $(0,\vep)$. In particular, there 
cannot be a first jump. Consequently some argument is needed
to justify the induction in the jump rule, so that 
formula (\ref{jump2}) can be used to define the successive
locations of particle $\sigma(i)$. 

 The realizations $\{\DD_i^k\}$ for which we can construct
the dynamics  satisfy  these assumptions: 
\bea
\mbox{(a)}&&\mbox{Each $\DD_i$ has only finitely many epochs 
in each bounded }\nn\\ 
 &&\mbox{time interval.}  \nn\\
\mbox{(b)}&&\mbox{There are no simultaneous jump attempts.}
\nn\\
\mbox{(c)}&&\mbox{There are arbitrarily large $i_0$ and $t_0$ such that 
$\DD_{-i_0}$ and $\DD_{i_0}$ }\nn\\
&&\mbox{have no epochs in $[0,t_0]$.}
\label{assD}
\eea
It is a standard fact of Poisson processes that 
these assumptions are satisfied by almost every realization 
$\{\DD_i^k\}$. 

Assumptions (\ref{assD})  allow us to prove that the jump rule
defines the evolution unambiguously for all
particles $\sigma(i)$ and all times $0\le t<\infty$:
 Given $i$ and $t$, part (c) of (\ref{assD}) gives 
$i_0$ and $t_0$ so that $-i_0<i<i_0$ and $t_0>t$, and 
so that $\DD_{-i_0}$ and $\DD_{i_0}$ 
have no epochs in $[0,t_0]$. Consequently 
 $\sigma(-i_0,t)=\sigma(-i_0,0)$ and 
$\sigma(i_0,t)=\sigma(i_0,0)$ for all $t\in[0,t_0]$,
and  the evolution
of particles $\sigma(j)$ for $-i_0<j<i_0$ is isolated 
from the rest of the process up to time $t_0$. By
parts (a) and (b) of (\ref{assD}), the 
locations  $\sigma(j,t)$ for $-i_0<j<i_0$ and  $t\in[0,t_0]$
can be computed by applying the jump rule 
to the finitely many potential jump times  
in $\cup_{-i_0<j<i_0}\DD_j\cap[0,t_0]$, in their temporal order. 
In particular, 
the motion of particle $\sigma(i)$ is defined up to time $t$. 

This way the evolution $\sigma(\cdot)=\{\sigma(t):t\ge 0\}$ is defined as a 
function of an arbitrary initial configuration $\sigma(0)$, and 
of almost every realization of the jump times $\{\DD^k_i\}$. 
If 
$(\Omega, {\cal F}, P)$ is  a probability space on which are defined the
 point processes 
$\{\DD^k_i\}$ and a random  initial configuration 
$\sigma(0)$,  then $\sigma(\cdot)$ is a 
stochastic process defined on this same probability space. 
Note: The Poisson point processes have to be independent 
of the initial  configuration 
$\sigma(0)$. 

Constructing particle systems from Poisson processes
of potential jump times  is a standard procedure, and 
 can be found in 
general references on particle systems such as 
 \cite{durrett, gr, lg}. 
A useful property of this construction is that it preserves
orderings between two processes:

\begin{lm} Suppose the probability space $(\Omega,{\cal F}, P)$
contains two initial configurations $\sigma'(\cdot\,,0)$
and $\sigma''(\cdot\,,0)$ that satisfy 
$\sigma'(i,0)\le \sigma''(i,0)$ for all $i\in\ZZ$ with probability $1$. 
If the processes  $\sigma'(\cdot\,,t)$
and $\sigma''(\cdot\,,t)$ obey the same point processes
$\{\DD^k_i\}$, we have 
$\sigma'(i,t)\le \sigma''(i,t)$ for all $i\in\ZZ$
and $t\ge 0$, with probability $1$. 
\label{orderinglm}
\end{lm}

{\it Proof.}
Choose $i_0$ and $t_0$
to satisfy part (c) of assumption (\ref{assD}). We shall prove 
that the ordering $\sigma'(i,t)\le \sigma''(i,t)$
holds for $-i_0\le i\le i_0$
and $t\le t_0$. It holds for $i=\pm i_0$ because 
particles $\sigma'(\pm i_0)$ and  $\sigma''(\pm i_0)$
do not attempt to jump during the time interval $[0,t_0]$.
 By
parts (a) and (b) of (\ref{assD}), we may do 
induction over the jump times. 
Suppose $\tau\in\DD^k_i$ is the first epoch at which the 
ordering is violated, for  $-i_0\le i\le i_0$
and $\tau\le t_0$. Then we must have $\sigma'(i,\tau)>\sigma''(i,\tau)$
but   $\sigma'(i,\tau-)\le\sigma''(i,\tau-)$.
Since  
$\sigma''(i)$ was not allowed to jump as far as 
$\sigma'(i)$ even though they both attempted a $k$-jump, 
it must be that   $\sigma'(i+1,\tau-)>\sigma''(i+1,\tau-)$.
Since $\tau$ is the first time the ordering
is violated for indices $-i_0\le i\le i_0$, it must be  that 
$i+1>i_0$. Since $i_0\ge i$, we conclude that $i=i_0$. But
 $\DD_{i_0}$ was chosen to have 
no epochs  up to time $t_0$. We have reached a 
contradiction, so there can be no first violation of the 
ordering for  $-i_0\le i\le i_0$
and $t\le t_0$. Since $i_0$ and $t_0$ can be taken 
arbitrarily large, the proof is complete. 
\hspace*{1mm} \rule{2mm}{3mm}

\section{The heuristic behind the method}
\label{heur}

We are interested in a macroscopic view of the process. By
this we mean its 
large scale behavior, distinct from the 
local interactions of the individual particles 
 $\sigma(i)$ that represent the microscopic
dynamics. 
 The ratio between macroscopic and microscopic scales
is given by a parameter $n$, so that macroscopic space and time units 
represent $n$ microscopic space and time units. 
We hope that a deterministic macroscopic
 description of the dynamics
emerges in the limit $n\to\infty$, where the
 macroscopic and microscopic scales become infinitely
separated. 

What form should the deterministic macroscopic 
description take? Suppose we have a function 
$u(x,t)$ defined for macroscopic space-time variables
$(x,t)$ such that the  approximate equality
\be
u(x,t)\approx n^{-1}\sigma([nx],nt)
\label{appr0}
\ee
 holds
with high probability for large $n$.
Since the  rate of 
 advance of  particle
$\sigma([nx],nt)$ depends only on the distance to the 
next particle on the right, we might expect this
same structure  also at the macroscopic
level, so that  
\be
u_t=f(u_x)
\label{hj}
\ee
for some function $f$. Equation (\ref{hj}) is a 
first-order partial differential equation of the  
Hamilton-Jacobi type. It is usually written as 
$u_t-f(u_x)=0$, and then the function $-f$ is known
as the {\it Hamiltonian}. 

 The exclusion rule 
and (\ref{appr0}) force $u_x\ge 1$, so $f$ is a function
defined on $[1,\infty)$. If the distance to the next
particle ahead
is  one (no empty sites) 
a particle cannot jump, and hence $f(1)=0$. The 
rate of advancement cannot decrease as the distance
to the next particle increases, so $f$ ought to be nondecreasing. 
No particle can advance faster than 
at rate $B_1$, so $f\le B_1$. Since $f$ is  nondecreasing
and bounded above, it appears reasonable to assume 
 that $f$ is concave.

Thus, assuming our heuristic reasoning is reliable,  if there
is to be a law of large numbers 
\be
n^{-1}\sigma([nx],nt)\to u(x,t)\,,
\label{lln}
\ee
 then $u$ should satisfy
equation  (\ref{hj}) with a convex Hamiltonian $-f$. 
There is a well-known formula for solving  (\ref{hj})
when the Hamiltonian is convex:
Let 
\be
g(x)=\sup_{ v\ge 1}\{ v x+f( v)\}
\label{gfdual}
\ee
be the convex conjugate of $-f$. Suppose the initial 
data for  (\ref{hj}) is given by a function $u_0$ on $\RR$. 
Then the {\it Hopf-Lax formula} 
\be
u(x,t)=\inf_{y\ge x}\left\{ u_0(y)+tg\left(\frac{x-y}t\right)\right\}
\label{hl}
\ee
defines the unique viscosity solution of  (\ref{hj}) with
initial data $u(x,0)=u_0(x)$ 
(Theorem 3, Section 10.3 in \cite{ev}). Note: Since $f\ge 0$, 
(\ref{gfdual}) shows that $g(x)=\infty$ for $x>0$. 
That is why we may restrict the infimum in (\ref{hl}) to
$y\ge x$. 

However, our task is not to solve 
the differential equation ({\ref{hj}) with a given
 $f$. In some sense we already have the solutions, at least in their 
microscopic
form $\sigma([nx],nt)$ because we already constructed the
particle process.
Instead, our problem is to prove that the 
process $\sigma$ admits a macrosopic description
of the type ({\ref{hj}), and to find $f$ and $g$, if possible. 
The Hopf-Lax formula
suggests a line of attack on our problem. 

Suppose we are allowed to feed a single
 input $u_0=\varphi_0$ to  the system
described by  (\ref{hl}),  and measure the output
$\varphi(x,t)=u(x,t)$. Can we obtain $g$? Take 
\be
\varphi_0(y)=\left\{ \begin{array}{lc}
\infty &\mbox{if $y>0$,}\\
y &\mbox{if $y\le 0.$}
\end{array}\right.
\label{specialu}
\ee
 With $u_0=\varphi_0$, 
differentiating inside the braces in (\ref{hl}) 
shows that $\varphi(x,t)=u(x,t)=tg(x/t)$. 
(The duality (\ref{gfdual}) implies that the slope of
$g$ is everywhere at least 1, since this is the lower bound
for $ v$.)
 Once we have 
$g$, we can compute the solutions to (\ref{hj})
by the Hopf-Lax formula (\ref{hl}). In other words, the
single evolution $\varphi(x,t)$ started with  (\ref{specialu}) contains 
 the information for constructing {\it all} evolutions
from their initial data. 

To understand this point better, let $u_0$ be
an arbitrary initial function for the evolution, 
and consider this family of 
special initial data, indexed by $z\in\RR$: 
\be
\varphi_0^z(y)= u_0(z)+\varphi_0(y)\,,\quad y\in\RR\,.
\label{specialu1}
\ee
These are simply  translates of the initial data (\ref{specialu}). 
Find the evolution from (\ref{hl}):
\beas
\varphi^z(x,t)&=&\inf_{y\ge x}
\left\{ \varphi^z_0(y)+tg\left(\frac{x-y}t\right)\right\}\\
&=& u_0(z)+\inf_{y\ge x} 
\left\{ \varphi_0(y)+tg\left(\frac{x-y}t\right)\right\}\\
&=&  u_0(z)+tg(x/t)\,.
\eeas
We substitute the evolution
$\varphi^z(x,t)$ into (\ref{hl}) to obtain
\be
u(x,t)=\inf_{y\ge x} \varphi^y(x-y,t)\,.
\label{hl1}
\ee
Equations (\ref{specialu1}) and (\ref{hl1}) contain the
message: A general evolution $u(x,t)$ starting from initial
data $u_0$ can be computed from the family of solutions
$\{\varphi^y:y\in\RR\}$, whose initial data are translates
 $\varphi_0^y= u_0(y)+\varphi_0$ of the 
function $\varphi_0$. 

Let us move the discussion to the microscopic particle level. 
We can realize the initial data  (\ref{specialu}) 
microscopically, in the sense of the limit (\ref{lln}), 
 by a natural particle
configuration:
\be
\sigma(i,0)=\left\{ \begin{array}{lc}
\infty &\mbox{if $i>0$,}\\
i &\mbox{if $i\le 0$.}
\end{array}\right.
\label{specialsigma}
\ee
This says that the lattice is 
empty to the right of the origin, while to the left all sites
are occupied, the origin itself included. 
Since  (\ref{specialsigma}) is a microscopic version 
of  (\ref{specialu}), 
 we would expect that the 
 process $\sigma(t)$ with initial configuration
 (\ref{specialsigma}) is a  microscopic version 
of the solution $\varphi(x,t)= tg(x/t)$.
Translates of these solutions  suffice to  express the most general 
solution, as observed  in (\ref{hl1}).

 The key question
becomes: Is there  a  microscopic  principle that
corresponds to  (\ref{hl1}) and 
expresses a process with  a general initial 
 configuration in terms of a family of 
processes of the type that 
 start from (\ref{specialsigma})? 
The construction we seek is the coupling we now turn to.

\section{The coupling}
\label{coupsec}

Assume given an arbitrary deterministic initial
particle configuration $\sigma(0)=(\sigma(j,0):j\in\ZZ)$,
and let $\sigma(j,t)$ denote the process constructed
in terms of the Poisson point processes $\{\DD^k_j\}$
as explained in Section \ref{constr}.
For each finite initial location $\sigma(j,0)$ we construct an
 auxiliary exclusion process 
 $\zeta^j=(\zeta^j(i,t):i\in\ZZ_-, t\ge 0)$. Initially
\be
\mbox{$\zeta^j(i,0)=\sigma(j,0)+i$ for $i\in\ZZ_-$.}
\label{zeta1}
\ee
The particles $\zeta^j(i)$ for $i>0$ are not needed. 
According to our convention 
we may think that they reside permanently at $\infty$. 
Thus $\zeta^j$ is a process that is a spatial translate of
the process started from (\ref{specialsigma}). 

For each
process $\zeta^j$, 
the dynamical description is the same as for $\sigma$: Jumps of size $k$
are attempted independently at rate $\beta_k$.  
Each jump is carried as far as possible without 
violating the exclusion and no-passing rules. Each
process $\zeta^j$ uses the same collection of Poisson processes
 $\{\DD^k_i\}$ as does $\sigma$, but with a translation of the 
index: 

\paragraph{Jump rule for process $\zeta^j$.}
At each epoch $\tau$ of $\DD^k_{i+j}$, 
 set
\be
\zeta^j(i,\tau)=\min\left\{ \zeta^j(i,\tau-)+k\,,\,
\zeta^j(i+1,\tau-)-1
\right\}\,.
\label{jumpzeta}
\ee
Between epochs of  $\DD_{i+j}$,  $\zeta^j(i)$
 stays constant. 

\hbox{}

 Subject to
this rule, the processes $\{\zeta^j\}$ are constructed
exactly as $\sigma$. The processes $\sigma$ and $\{\zeta^j\}$ are
all defined on the same probability space
$(\Omega,{\cal F},P)$ of 
the  Poisson jump times. If this probability space also
supports a random initial configuration $\sigma(0)$, then 
the processes $\{\zeta^j\}$ are functions of both 
 $\sigma(0)$ and  $\{\DD^k_i\}$, through the translations
(\ref{zeta1}) and the jump rule.

 This arrangement where
several processes are defined on a common probability space
to facilitate their direct comparison 
is known as a {\it coupling}. 
 It should be evident that
the processes are invisible to each other, i.e.\ particles of
one process do not interfere with the jumps of the other processes. 

The key property of this coupling is the microscopic 
counterpart of (\ref{hl1}): 

\begin{lm} The equality
\be\sigma(i,t)=\inf_{j:j\ge i}\zeta^{j}(i-j,t)
\label{coup1}
\ee
holds for all $i\in\ZZ$ and  $t\ge 0$, almost
surely. 
\label{couplinglm}
\end{lm}

{\it Proof.} The exclusion rule (\ref{e1}) and 
(\ref{zeta1}) imply 
 that (\ref{coup1})
holds at time 0.  The jump rules
ensure that for each $j\ge i$, particles $\zeta^j(i-j)$
and $\sigma(i)$ attempt $k$-jumps at common time points,
namely at the epochs of $\DD^k_i$. Thus
for a fixed $i$ the validity of 
(\ref{coup1}) can change only at epochs of $\DD_i$. 

Fix $i_1$ and $t_1$. To prove that (\ref{coup1}) holds
for $i=i_1$ up to time $t_1$, use part (c) of assumption 
(\ref{assD}) to find $i_0$ such that $-i_0<i_1<i_0$, and so that
$\DD_{-i_0}$ and $\DD_{i_0}$ have no epochs in
the time interval $[0,t_1]$. By  parts (a) and (b) of (\ref{assD}),
we can do induction on the finitely many epochs in
$\cup_{-i_0<i<i_0}\DD_i\cap[0,t_1]$. So suppose 
$\tau$ is an epoch of $\DD^k_i$ for some $k\ge 1$ 
and some $i$ between $-i_0$ and $i_0$, and  assume 
that  (\ref{coup1}) holds for all  $i$ between $-i_0$ and $i_0$
and all $t<\tau$. 

First we show that $\le$ holds in  (\ref{coup1}) right after
the jump epoch $\tau$. Let $j\ge i$ be arbitrary. 
By the induction assumption, $\zeta^j(i-j+1,\tau-)\ge 
\sigma(i+1,\tau-)$, and consequently $\zeta^{j}(i-j)$ can
 jump at least as far as $\sigma(i)$ can. Since
before the jump we have  $\zeta^j(i-j,\tau-)\ge 
\sigma(i,\tau-)$, and both particles attempt a $k$-jump,
the inequality  $\zeta^j(i-j,\tau)\ge 
\sigma(i,\tau)$ holds also at time $\tau$  after the 
jump has been completed. The same argument works for all
 $j\ge i$. 

Next we show that right after $\tau$ there is some
$j\ge i$ that satisfies  $\zeta^j(i-j,\tau)=  
\sigma(i,\tau)$. If $\sigma(i)$ was allowed to
jump the full $k$ steps, any $j$ such that 
 $\zeta^j(i-j,\tau-)=  
\sigma(i,\tau-)$ will do, and by induction there is some
such $j$. 

Suppose on the contrary that $\sigma(i)$ was
at least partially blocked by $\sigma(i+1)$ at time $\tau-$, so after
the jump we have  $\sigma(i,\tau)=\sigma(i+1,\tau)-1$.
By induction, we can pick  $j\ge i+1$ such that
 $\zeta^j(i-j+1,\tau-)=  
\sigma(i+1,\tau-)$. By the earlier part
of the proof,  $\zeta^j(i-j,\tau)\ge 
\sigma(i,\tau)$. On the other hand, the exclusion
rule for the process $\zeta^j$ stipulates that 
 $\zeta^j(i-j,\tau)$ $\le$ $\zeta^j(i-j+1,\tau)-1$
$=$ $\sigma(i+1,\tau)-1$ $=$ $\sigma(i,\tau)$. 
(Notice that by part (b) of assumption (\ref{assD}) there is 
no epoch $\tau$ in $\DD_{i+1}$.)
 This shows that after the jump we necessarily have the 
equality   $\zeta^j(i-j,\tau)=  
\sigma(i,\tau)$. 
\hspace*{1mm} \rule{2mm}{3mm}

\vspace*{3mm}

By subtracting off the centering at $\sigma(j,0)$
 from $\zeta^j$ we obtain
a family of processes $\{\xi^j\}$ that are all equal in distribution.
Define 
\be
\mbox{$\xi^j(i,t)=\zeta^j(i,t)-\sigma(j,0)$ for all 
$i,j$, and $t$.}
\label{defxi}
\ee
Then we can turn the coupling equality (\ref{coup1})
into a statement that is a microscopic version of the Hopf-Lax
formula (\ref{hl}):
\be
\sigma(i,t)=\inf_{j:j\ge i}\left\{ \sigma(j,0)+\xi^{j}(i-j,t)\right\}.
\label{coup2}
\ee
 The processes $\{\xi^j\}$ are identically distributed on
account of the requirement that $\left(\sigma(j,0)\right)$ and
the Poisson processes $\{\DD^k_j\}$ are independent. 
The advantage of (\ref{coup2}) over  (\ref{coup1}) 
is that the  effects of  $\left(\sigma(j,0)\right)$ and  $\{\DD^k_j\}$
 have been  separated in two distinct terms
inside the braces.

\section{The tagged particle limit}
\label{limit}

As a  corollary of the coupling we obtain the limit discussed
heuristically in Section \ref{heur}. The precise 
assumptions are as follows: There is a sequence of 
processes $\sigma_n$ constructed as in Section \ref{constr},
and an increasing function $u_0$ on $\RR$. 
At time 0, the  convergence  
\begin{equation}
\lim_{n\to\infty} n^{-1}\sigma_n([ny],0)=u_0(y)
\label{lim1}
\end{equation}
holds in probability for each
$y\in\RR$. More precisely, if $(\Omega_n, \FF_n, P_n)$ is the 
probability space on which the process $\sigma_n$ is defined,
then 
\begin{equation}
\lim_{n\to\infty} P_n\left( \left|\sigma_n([ny],0)-nu_0(y)
\right|\ge n\vep\right)=0
\label{lim1a}
\end{equation}
for each
$y\in\RR$ and $\vep>0$.

\begin{thm} Under assumption {\rm (\ref{lim1})},
 the convergence  
\begin{equation}
\lim_{n\to\infty} n^{-1}\sigma_n([nx],nt)=u(x,t)
\label{lim2}
\end{equation}
holds in probability for each $x\in\RR$ and $t>0$. 
The deterministic limit $u(x,t)$ can be 
macroscopically defined by 
\be
u(x,t)=\inf_{y:y\ge x}\left\{ u_0(y)+tg\left(\frac{x-y}{t}\right)
\right\},
\label{defU}
\ee
where $g$ is a certain convex function on $(-\infty,0]$ 
determined by the rules of the exclusion process.
\label{thm1}
\end{thm}

\vspace*{0.3cm}

The meaning of the limit in (\ref{lim2}) is the same as 
above for  (\ref{lim1}), namely that 
\begin{equation}
\lim_{n\to\infty} P_n\left( \left|\sigma_n([nx],nt)-nu(x,t)
\right|\ge n\vep\right)=0
\label{lim2a}
\end{equation}
for each $\vep>0$.
By the discussion about Hamilton-Jacobi equations in
Section \ref{heur}, $u(x,t)$ can be equivalently characterized
as the unique viscosity solution of $u_t=f(u_x)$ with
initial data $u(x,0)=u_0(x)$, where $f$ is the negative
of the convex conjugate of $g$. The theorem says that
in this exclusion system, a tagged particle has a 
well-defined macroscopic velocity $f(u_x)$ as long as 
initially the particle distribution is sufficiently
regular as specified by assumption (\ref{lim1}). 

\hbox{}

\subsubsection*{Proof of Theorem 1}
To prove Theorem \ref{thm1},   we forge a rigorous 
connection between the matching microscopic
and macroscopic descriptions (\ref{coup2}) and (\ref{defU}). 
For each $n$, we have the processes $\zeta_n^j$ and $\xi_n^j$,
defined as in Section \ref{coupsec} relative to the process
$\sigma_n$. All the processes $\xi_n^j$ are identical 
in distribution, so as long as we discuss distributional 
properties, we may leave out the sub- and superscript $n$ and  $j$. 
The function $g$ that appears in the theorem is defined by
the limit of the next lemma. 

\begin{lm} There exists a finite, continuous, 
 convex function $g$ on $(-\infty,0]$ such that 
\begin{equation}
\lim_{n\to\infty} n^{-1}\xi([nx],nt)=tg(x/t)
\label{limxi}
\end{equation}
 in probability,  for each $x\le 0$ and $t>0$. 
\label{xilemma}
\end{lm}

{\it Proof of Lemma \ref{xilemma}.} 
We start by deriving  a useful subadditivity.
Assume the process $\xi(i,t)$ has been constructed 
on a probability space $(\Omega,\FF,P)$ 
as explained in Section \ref{constr}, so that particle 
$\xi(i)$ gets its jump commands from the labeled Poisson process
$\DD_i$. To perform comparisons,  we explicitly 
include the particles that are permanently at infinity:
$\xi(i,t)=\infty$ for all $i>0$ and $t\ge 0$. 

 Fix  an integer $h\in\ZZ_-$  and a time $s>0$. Construct two new processes 
$\sigma'$ and $\sigma''$
on the same  probability space $(\Omega,\FF,P)$
 as follows: The initial configurations
are 
\be 
\mbox{$\sigma'(i,0)=\xi(h+i,s)$ for all $i\in\ZZ$}
\label{sigmap0}
\ee
and 
\be 
\sigma''(i,0)=\left\{ \begin{array}{lc}
\infty &\mbox{if $i>0$,}\\
\xi(h,s)+i &\mbox{if $i\le 0$.}
\end{array}\right.
\label{sigmapp0}
\ee
We stipulate that particles
 $\sigma'(i)$ and $\sigma''(i)$ read their jump commands
from $\left(\DD_{h+i}-s\right)\cap[0,\infty)$.  In other words,
$\sigma'$ and $\sigma''$ read
 the Poisson processes $\{\DD_i\}$ after first translating
the index by $h$ and time by $s$. By Lemma
\ref{orderinglm}, 
\be
\mbox{$\sigma'(i,t)\le \sigma''(i,t)$ for all $i\in\ZZ$ and $t\ge 0$}.
\label{ordersigma}
\ee
Since $\sigma'$ continues $\xi$ with translated 
index and time, $\sigma'(i,t)=\xi(h+i,s+t)$. For
 $\sigma''$ we can write  
$\sigma''(i,t)=\xi(h,s)+\xit(i,t)$ where $\xit$ is equal
in distribution to $\xi$ but independent of it. Thus
(\ref{ordersigma}) turns into
\be
\xi(h+i,s+t) \le\xi(h,s)+\xit(i,t) \,,
\label{subadd}
\ee
where the random variables on the right are independent.

Suppose first that $x$ in (\ref{limxi})
 is a nonpositive integer, $x=i$. Let $F_n$ be the 
distribution function of $\xi(ni,nt)$. Then 
 (\ref{subadd}) gives $F_{n+m}\ge F_n*F_m$. By 
the Kesten-Hammersley lemma from subadditive ergodic 
theory there exists a 
function $\gamma(i,t)$ such that
\begin{equation}
\lim_{n\to\infty} n^{-1}\xi(ni,nt)=\gamma(i,t)
\label{limxi1}
\end{equation}
 in probability,  for all $i\in\ZZ_-$ and $t\ge 0$. 
A proof of this lemma can be found on p.\ 20 in
\cite{sw}. The assumption of finite second moments
required by this lemma 
is satisfied 
 by assumption (\ref{assB2}). 

The rest of the proof consists of extending the limit
to all $x\le 0$, by using regularity and monotonicity
properties of the process and $\gamma(i,t)$. 
It is obvious that $\gamma(i,t)$ is nondecreasing in 
both $i$ and $t$. For $i=0$ we have 
$\gamma(0,t)=B_1t$, because $B_1$ is the speed of an 
unobstructed particle. Furthermore 
\be
i\le \gamma(i,t)\le i+B_1t
\label{boundgamma}
\ee
 for all $i\in\ZZ_-$ and $t\ge 0$, because 
$\xi(i,0)=i$ and $B_1$ is the maximal average speed. 

Next we argue the $t$-continuity of $\gamma(i,t)$. 
During the time interval $(nt,nt+n\vep]$
 particle   $\xi(ni)$
attempts a Poisson($B_0n\vep$) number of jumps. The 
attempted jump sizes are i.i.d.\ random variables  $\{X_i\}$ 
with common  distribution $P(X_i=k)=\beta_k/B_0$. Due to the 
exclusion and no-passing rules, the actual jumps are
at most $X_i$, so we can write 
\be
\xi(ni, nt+n\vep)\le \xi(ni,nt) +\sum_{i=1}^{N}X_i
\label{cont1}
\ee
where $N\sim$ Poisson($B_0n\vep$). 
 It follows that
\be
\gamma(i,t+\vep)\le \gamma(i,t)+B_1\vep\,.
\label{contgamma}
\ee
In particular, $\gamma(i,t)$ is 
Lipschitz continuous in $t$. 

A homogeneity holds: If $m$ is a positive integer, then
\beas
\gamma(mi,mt)&=&\lim_{n\to\infty} n^{-1}\xi(nmi,nmt) \\
&=&m\cdot\lim_{n\to\infty} (nm)^{-1}\xi(nmi,nmt) \\
&=&m\gamma(i,t)\,.
\eeas
From this it follows that we can define unambiguously 
for rational $r\le 0$ 
\be
\gamma(r,t)=\frac1m \gamma(mr,mt)\,,
\label{defrat}
\ee
where $m$ is any positive integer such that $mr\in\ZZ_-$. 

We need  to prove that the extension (\ref{defrat}) 
makes sense for the process, namely that 
\be
\lim_{n\to\infty} n^{-1}\xi([nr],nt)=\gamma(r,t)
\label{limrat}
\ee
 for rational $r< 0$
and all $t>0$. If $n$ runs along the subsequence 
$n=km$, $k=1,2,3,\ldots$ , where $m$ is the integer
in  (\ref{defrat}), then (\ref{limrat}) is valid. 
For other values of $n$ we interpolate using the monotonicity
of the process. No-passing and total asymmetry
of the dynamics imply that 
$$\xi(i,s)\le \xi(j,s)\le \xi(j,t)$$
for all $i\le j$ and $s\le t$. 
Given $n$, let $k=k(n)$ be the unique 
integer that satisfies $km\le n<(k+1)m$. Let $\vep>0$. 
Then for large enough $n$,
\beas
\xi\left((k+1)mr,(k+1)m(t-\vep)\right)&\le&\xi([nr],nt)\\
&\le&\xi\left(kmr,km(t+\vep)\right)\,,
\eeas
from which, upon dividing by $n$ and passing to the limit, 
\beas
\frac1m \gamma\left(mr, m(t-\vep)\right)&\le&
\liminf_{n\to\infty} n^{-1}\xi([nr],nt)\\
&\le&
\limsup_{n\to\infty} n^{-1}\xi([nr],nt)\\
&\le& \frac1m \gamma\left(mr,m(t+\vep)\right)\,.
\eeas
These inequalities cannot be taken literally
since the limits  have to be interpreted in the sense
of convergence in probability. But the point should
be clear: By the continuity (\ref{contgamma}) and the
definition (\ref{defrat}), the limit (\ref{limrat}) is valid.

Next one checks that the homogeneity extends to rationals:
\be
\gamma(rx,rt)=r\gamma(x,t)
\label{homogamma}
\ee
for rational $x\le 0$, rational $r\ge 0$, and all $t>0$. 
From (\ref{subadd}) follows  subadditivity
\be
\gamma(x+y,s+t)\le \gamma(x,s)+\gamma(y,t)\,,
\label{subaddgamma}
\ee
first for integers $x,y$ and then by  (\ref{defrat})
for rational  $x,y$. (\ref{homogamma}) and (\ref{subaddgamma})
together imply convexity, and also continuity in the 
$x$-variable, for rational $x$. 

The final extension of $\gamma$ is
\be
\gamma(x,t)=\sup\{\gamma(r,t): r\in\QQ, r<x\}\,.
\label{defirr}
\ee
By the continuity in the $x$-variable, 
for rational $x$ this definition gives the old value
$\gamma(x,t)$, so it is a sensible extension. 
Again, one checks that homogeneity (\ref{homogamma})
and  subadditivity (\ref{subaddgamma}) 
hold, this time for all $r,x,y$. Convexity follows,
and  convexity implies
continuity in the open quadrant $\{x<0, t>0\}$
(Theorem 10.1 in \cite{ro}).
Continuity up to the boundary requires separate 
arguments. And again, one  checks that the limit 
(\ref{limrat}) holds also for irrational $r$. This 
follows as it did for rational $r$, by the 
monotonicity of the process and the continuity
of the limit. 

Finally, define $g(x)=\gamma(x,1)$. By homogeneity
$\gamma(x,t)=tg(x/t)$. 
\hspace*{1mm} \rule{2mm}{3mm}

\vspace*{3mm}

We return to the proof of Theorem \ref{thm1}. 
Fix $(x,t)$. From the coupling (\ref{coup2}) we get 
\be
n^{-1}\sigma_n([nx],nt)=
\inf_{i:i\ge [nx]}\left\{n^{-1} \sigma_n(i,0)
+n^{-1}\xi_n^{i}\left([nx]-i,nt\right)\right\}.
\label{coup3}
\ee
By the hypothesis of Theorem \ref{thm1} and by
Lemma \ref{xilemma}, for $i=[ny]$ 
the random variable inside the braces converges
to the quantity inside the braces in (\ref{defU}). 
Consequently we get 
\be
\lim_{n\to\infty}
P_n\left(n^{-1}\sigma_n([nx],nt)\le u(x,t)+\vep\right)=1
\label{thm11}
\ee
for any $\vep>0$. 

For the converse we need some further  estimation. 
The first step is to restrict the variable $i$
in (\ref{coup3}) to a range of order $n$. This is achieved
by the next lemma. Define 
\be
w_n(z)=\min_{i:[nx]\le i\le [nz]}\left\{ \sigma_n(i,0)
+\xi_n^{i}\left([nx]-i,nt\right)\right\}.
\label{defw}
\ee

\begin{lm} For large enough $z>x$, there exists
a constant $C=C(z)>0$ such that 
\begin{equation}
P_n\left(\sigma_n([nx],nt)\ne w_n(z)\right)\le e^{-Cn}
\label{wbound}
\end{equation}
for all large enough $n$. 
\label{wlemma}
\end{lm}

{\it Proof.} Suppose particle 
$\xi_n^{[nz]}\left([nx]-[nz]\right)$ has not moved by time 
$nt$, in other words that  $\xi_n^{[nz]}\left([nx]-[nz],nt\right)$
$=$ $[nx]-[nz]$. Let $i>[nz]$. 
By the exclusion rule for $\sigma_n$  and because 
$$\xi_n^{i}\left([nx]-i,nt\right)\ge
\xi_n^{i}\left([nx]-i,0\right) =[nx]-i\,,
$$ 
we have 
\beas
&& \sigma_n(i,0)
+\xi_n^{i}\left([nx]-i,nt\right)\\
&\ge& \sigma_n\left([nz],0\right)+i-[nz]
+[nx]-i\\
&=& \sigma_n\left([nz],0\right)+[nx]-[nz]\\
&=& \sigma_n\left([nz],0\right)+\xi_n^{[nz]}\left([nx]-[nz],nt\right)\,.
\eeas
The conclusion is that the terms  $i>[nz]$
cannot contribute to the infimum in (\ref{coup3}). 

Thus to prove the lemma, it suffices to show that
for large enough $z$,
the probability that $\xi_n^{[nz]}\left([nx]-[nz]\right)$ 
has moved by time 
$nt$ is at most $e^{-Cn}$. This is clear again on account
of the exclusion and no-passing rules: The random time when 
 $\xi_n^{[nz]}\left([nx]-[nz]\right)$ first moves has
the distribution of a sum of $[nz]-[nx]+1$ i.i.d.\ exponential
random variables with rate $B_0$, because after
$\xi(j)$ first moves, $\xi(j-1)$ waits an Exp($B_0$)-distributed
independent time to make its first jump, and these times are added up
over $j=0,-1,-2$, $\ldots$, $[nx]-[nz]$. It is sufficient
to choose $z$ so that  $(z-x)B_0^{-1}>t+\vep$ for some
$\vep>0$. 
\hspace*{1mm} \rule{2mm}{3mm}

\vspace*{3mm}

Fix $z$ so that (\ref{wbound}) holds, and let $\vep>0$. 
Pick a partition 
$$x=y_0<y_1<y_2<\cdots<y_m=z$$
so that
\be
\left| tg\left(\frac{x-y_{k+1}}t\right) - 
 tg\left(\frac{x-y_{k}}t\right) \right|\le \frac{\vep}4
\label{part1}
\ee
for $k=0,1,\ldots,m-1$. If we could restrict the minimum
in (\ref{defw}) to the $m+1$ points
$i=[ny_k]$,  then passing to the limit inside or outside
the minimum (now over a fixed finite number of random
variables) would make no difference. To control
the error of this 
simplification, we need to control what happens inside
the braces in (\ref{defw}) for $[ny_k]< i< [ny_{k+1}]$. 
This can be done with couplings: 

\begin{lm} For all $n$ and $i\le j_0\le j_1$, 
the inequality $\xi_n^{j_0}(i-j_0,t)$ $\ge$ 
$\xi_n^{j_1}(i-j_1,t)$ holds almost surely for all times $t\ge 0$. 
\label{xilemma1}
\end{lm}

{\it Proof.} Fix $n$ and $j_0\le  j_1$. Consider the processes
$$
\mbox{$\sigma'(i,t)=\xi_n^{j_1}(i-j_1,t)$ and 
$\sigma''(i,t)=\xi_n^{j_0}(i-j_0,t)$, $i\in\ZZ$.}
$$
 Both 
$\sigma'(i)$ and $\sigma''(i)$ read their jump commands
from $\DD_i$. (Recall the general rule from Section \ref{coupsec}: 
particle $\xi^j(i)$ reads jump commands from $\DD_{i+j}$.) Initially
\be 
\sigma'(i,0)=\left\{ \begin{array}{lc}
\infty &\mbox{if $i>j_1$,}\\
i-j_1 &\mbox{if $i\le j_1$,}
\end{array}\right.
\label{sigmap1}
\ee
and 
\be 
\sigma''(i,0)=\left\{ \begin{array}{lc}
\infty &\mbox{if $i>j_0$,}\\
i-j_0 &\mbox{if $i\le j_0$.}
\end{array}\right.
\label{sigmapp1}
\ee
By the assumption  $j_0\le  j_1$, 
$\sigma'(i,0)$ $\le$ $\sigma''(i,0)$, and consequently
by Lemma  \ref{orderinglm}, the inequality 
$\sigma'(i,t)$ $\le$ $\sigma''(i,t)$ holds for all times $t\ge 0$.
\hspace*{1mm} \rule{2mm}{3mm}

\vspace*{3mm}

Apply Lemma \ref{xilemma1} with $i=[nx]$, $j_0=i$ between
$[ny_k]$ and $[ny_{k+1}]$, and $j_1=[ny_{k+1}]$, to
bound $w_n(z)$ from below by
\beas
w_n(z)&=&\min_{0\le k<m}
\,\min_{i:[ny_k]\le i\le [ny_{k+1}]}\left\{ \sigma_n(i,0)
+\xi_n^{i}\left([nx]-i,nt\right)\right\}\\
&\ge&\min_{0\le k<m}\left\{ \sigma_n\left([ny_k],0\right)
+\xi_n^{[ny_{k+1}]}\left([nx]-[ny_{k+1}],nt\right)\right\}\,.
\eeas
Multiplied through by $n^{-1}$, the last line
 converges 
in probability by assumption (\ref{lim1}) and
Lemma \ref{xilemma}. 
Thus with probability tending to $1$ as $n\to\infty$, 
\beas
n^{-1}w_n(z)&\ge&\min_{0\le k<m}\left\{
u_0(y_k)+ tg\left(\frac{x-y_{k+1}}t\right)\right\}-\frac{\vep}4\\
&\ge&\min_{0\le k<m}\left\{
u_0(y_k)+ tg\left(\frac{x-y_{k}}t\right)\right\}-\frac{\vep}2\\
&\ge&\inf_{y:y\ge x}\left\{
u_0(y)+ tg\left(\frac{x-y}t\right)\right\}-\frac{\vep}2\\
&=& u(x,t)-\frac{\vep}2\,.
\eeas
This together with Lemma \ref{wlemma} gives
\be
\lim_{n\to\infty}
P_n\left(n^{-1}\sigma([nx],nt)\ge u(x,t)-\vep\right)=1
\label{thm12}
\ee
for any $\vep>0$. 
This completes the proof of Theorem \ref{thm1}.

\section{Comments and extensions}
\label{comments}

\subsection{Almost sure convergence in Theorem \ref{thm1}}

It may be possible to upgrade
the  convergence in Theorem \ref{thm1} 
 to almost sure convergence by developing 
 large deviation estimates for the convergence
in Lemma \ref{xilemma}. Such results can be 
found in \cite{se2} and \cite{sk}.

\subsection{The role of invariant distributions}

Theorem 1 does not explicitly identify the function
$g$ or the macroscopic
velocity $f$. This can be done when information about
steady states is available. Consider the special case
of the totally asymmetric simple exclusion process
(TASEP) where $\beta_1=1$ and $\beta_k=0$ for $k>1$. 
Then it is known that if the interparticle distances
are geometrically distributed with mean $ v\in [1,\infty)$, 
\be
P\left(\sigma(i+1)-\sigma(i)=m\right)=\left(1-v^{-1}\right)^{m-1}
v^{-1}\,,\quad m=1,2,3,\ldots\,,
\label{eq1}
\ee
an individual particle's jumps obey a Poisson process
of rate $1-v^{-1}$. From this information the function
$g$ can be calculated with the help of Theorem \ref{thm1},
as we next show. 

Let each $\sigma_n=\sigma$ be a steady-state process
with interparticle distances (\ref{eq1}), and  such that
initially $\sigma(0,0)=0$. Then  Theorem \ref{thm1}
is valid for $u_0(y)=vy$.  Take $x=0$ and $t=1$
in (\ref{defU}) to get
\be
u(0,1)=\inf_{y\ge 0}\{ vy+g(-y)\}=-g^*(v)\,,
\label{gstar1}
\ee
where $g^*$ denotes the convex conjugate of $g$. 
Note: For the purposes of convex duality, functions should
be lower semicontinuous and convex on all of $\RR$.
(See \cite{ro}.) The limit of Lemma \ref{xilemma}
defines $g$ on $(-\infty,0]$, and we extend
$g$ to $\RR$ by defining $g(x)=\infty$ for $x>0$. Then 
(\ref{gstar1}) is the same as 
\be
g^*(v)=\sup_x \{xv-g(x)\}.
\label{dualg}
\ee
 Furthermore, for
$v<1$ we obtain $g^*(v)=\infty$ directly from (\ref{dualg}), 
by the observation that the slope
of $g$ is always at least 1. This is a consequence
of the exclusion rule (\ref{e1}) and the limit (\ref{limxi}). 

Since $\sigma(0,n)$ is a 
Poisson($(1-v^{-1})n$)-variable, we can calculate explicitly that
$$u(0,1)=\lim_{n\to\infty}n^{-1}\sigma(0,n)=1-v^{-1}.$$
Thus $g^*(v)=v^{-1}-1$ for $v\ge 1$.
 From this we obtain $g$ by performing
another convex conjugation: For $x\le 0$, 
\beas
g(x)&=&g^{**}(x)\\
&=&\sup_{v}\{xv-g^*(v)\}\\
&=&\sup_{v\ge 1}\{xv-v^{-1}+1\}\\
&=&1-2\sqrt{-x}\,.
\eeas

The tagged particle velocity is the function $f(v)=1-v^{-1}$
obtained above. The relationship of the mean
interparticle  distance $v$ to
the particle density $\rho$ is $\rho=v^{-1}$. 
 As a function of density the velocity is $v(\rho)=f(\rho^{-1})
=1-\rho$. From this  we get the well-known
macroscopic {\it particle current} of TASEP: 
$\rho v(\rho)=\rho(1-\rho)$. 

TASEP was the process studied by Rost in 1981 \cite{rs}
in one of the seminal papers on hydrodynamic limits for
asymmetric processes. With an approach completely
different from ours he  investigated the 
process $\xi$ and  proved a limit theorem related 
to  Lemma \ref{xilemma}. The use  of subadditivity is
a common feature of both approaches. 

A thorough treatment of TASEP with the 
method of this paper appears in \cite{se1}.
In particular, when the invariant distributions are
available for explicit calculations, in addition
to the law of large numbers it is possible to
obtain a (one-sided) explicit large deviation rate function 
for the tagged particle in a nonequilibrium process.  

\subsection{More general jump rates}

We could admit the following additional feature 
to the process: Not only
is  $\sigma(i)$ blocked by the particle in front of him,
but he is also held back by the particle behind
him, in this particular sense:  
For some constant $L\ge 1$, 
\begin{equation}
 \sigma(i)\le \sigma(i-1)+L \mbox{ for all $i\in\ZZ$.}
\label{L2}
\end{equation}
The proof of Theorem \ref{thm1} 
works with suitable modifications. The initial
 configurations of the $\zeta^j$-processes 
have to be changed to $\zeta^j(i)=\sigma(j,0)+i$
for $i\le 0$, and  $\zeta^j(i)=\sigma(j,0)+Li$
for $i> 0$. This feature appears in \cite{se2}. 

Another generalization that this approach can handle
is {\it random rates}. The reader may have noticed 
that very little depended on how the point processes
$\DD^k_i$ were actually produced, as long as they 
had some desirable properties such as those listed
in assumption (\ref{assD}). In particular, the rates
of the processes could vary so that the rate $\beta_{k,i}$
of $\DD^k_i$ is a random variable whose realization is
fixed. If this randomness is suitably systematic
(i.i.d.\ for example) one would expect that
some form of the subadditive ergodicity
needed for Lemma \ref{xilemma}  works again.
Examples appear in \cite{se2} and \cite{sk}. 

However, so far it has not been possible to admit a 
more complicated dependence of the rate of jumping 
on the distance to the nearest particles. By a more
complicated dependence we mean 
 a general function
$\beta(k,l,r)$ such that jumps of size $k$ are 
attempted at rate $\beta(k,l,r)$ when there are 
$l$ and $r$ empty sites to the left and right. 
In the model of this paper this function is a step function:
For a fixed $k$,  $\beta(k,l,r)$ takes only two values
depending on whether $k<r$ or $k\ge r$. 
What appears to fail for more general 
rates  is the proof of Lemma
\ref{couplinglm}. 

Other shortcomings are that the method
has been able to treat  only totally asymmetric
processes, and processes without passing.
 A more general asymmetric process
would allow jumps both left and right, but with 
a higher rate for jumping to the right.


\begin{thebibliography}{99}

\bibitem{durrett} R. Durrett (1995).
Ten  lectures on particle systems. Lecture Notes
in Mathematics 1608 (Saint-Flour, 1993), 97--201. Springer-Verlag. 

\bibitem{ev} L. C. Evans (1998).
Partial Differential Equations. American Mathematical Society.

\bibitem{fe} P. A. Ferrari (1996). Limit theorems for
tagged particles. Markov Process. Related Fields 2 17--40. 

\bibitem{gr}
D. Griffeath (1979). Additive
and Cancellative Interacting Particle Systems.
Lecture Notes in Mathematics 724, Springer, Heidelberg.


\bibitem{lg}
T. M. Liggett (1985). Interacting Particle Systems.
Springer-Verlag, New York.



\bibitem{ro} R. T. Rockafellar (1970).
Convex Analysis. Princeton University Press.

\bibitem{rs}
 H. Rost (1981).  Non-equilibrium behaviour of a many particle
process: Density profile and local equilibrium.
 Z. Wahrsch. Verw. Gebiete  58, 41--53.



\bibitem{se1}
T. Sepp\"al\"ainen (1998). Coupling the totally asymmetric
simple exclusion process with a
moving interface. Markov Process.\ Related Fields 4, 593--628.
Proceedings of I Escola Brasileira de Probabilidade, IMPA, 1997. 

\bibitem{se2}
T. Sepp\"al\"ainen (1999). Existence of 
hydrodynamics for  the totally asymmetric
simple $K$-exclusion process.  Ann.\ Probab.\ 27, 361--415.

\bibitem{sk}
T. Sepp\"al\"ainen and J. Krug (1999). 
Hydrodynamics and platoon formation  for a totally
asymmetric exclusion model with particlewise disorder. 
J. Statist. Phys. 95, 529--571.  

\bibitem{sw}
R. T. Smythe and J. C. Wierman (1978). First-passage percolation
on the square lattice. Lecture Notes in Mathematics 671, Springer. 
\end{thebibliography}
\end{document}